\let\bd=\bf  % used for definitions
\let\a=\alpha
\let\b=\beta
\let\s=\sigma
\def\Z{{\hbox{\bf Z}}}
\def\F{{\hbox{\bf F}}}
\def\R{{\hbox{\bf R}}}
\def\Q{{\hbox{\bf Q}}}
\def\Hom{{{\rm Hom}}}
\def\End{{{\rm End}}}
\def\cj#1{{{\overline{#1}}}}

% For making diagrams
\def\diagram#1{{\normallineskip=8pt \normalbaselineskip=0pt
\matrix{#1}}}
\def\harr#1#2{\smash{\mathop{\hbox to .5in{\rightarrowfill}}
\limits^{\scriptstyle#1}_{\scriptstyle#2}}}
\def\varr#1#2{\llap{$\scriptstyle#1$}\left\downarrow
\vcenter to .5in{}\right.\rlap{$\scriptstyle#2$}}

\input plncs.cmm
\contribution{Formal Groups, Elliptic Curves, and Some Theorems of Couveignes}
\author {Antonia W.\ Bluher}
\address{National Security Agency, 9800 Savage Road, 
Fort George G.\ Meade, MD  20755-6000}
\abstract{
The formal group law of an elliptic curve has seen recent applications
to computational algebraic geometry in the work of Couveignes
to compute the order of an elliptic curve over finite fields of small
characteristic ([2], [6]).  The purpose of this paper is 
to explain in an elementary way how to associate a formal group law
to an elliptic curve and to expand on some theorems of Couveignes.
In addition, the paper serves as background for [1].
We treat curves defined over arbitrary fields, including
fields of characteristic two or three.  
The author wishes to thank Al Laing for a
careful reading of an earlier version of the manuscript and for many useful
suggestions.}

\titlea{1}{ Definition and construction of formal group laws}

Let $R$ be a commutative ring with a multiplicative identity 1
and let $R[[X]]$ denote the ring of formal power series of $R$.
In general it is not possible to compose two power series in
a meaningful way.  For example, if we tried to form the composition
$f\circ g$ with $f=1+\tau+\tau^2+\tau^3+\cdots$ and $g=1+\tau$
we would get
$$f\circ g=1+(1+\tau)+(1+\tau)^2+(1+\tau)^3+\cdots$$
The constant term is $1+1+1+\cdots$, which makes no sense.
But there are some cases where $f\circ g$ does make sense,
namely when $f$ is a polynomial {\it or} when the constant term of
$g$ is zero.  Let
$R[[X,Y]]=R[[X]][[Y]]$, the ring of formal power series in
two variables.  If $F\in R[[X,Y]]$ and $g,h\in \tau R[[\tau]]$
then
$$F(g,h)\quad{\rm makes\ sense\ and\ belongs\ to\ } R[[\tau]].$$
If in addition $F$ has a zero constant term, then $F(g,h)\in\tau
R[[\tau]]$.
  
A one dimensional (commutative)
{\bd formal group law} over $R$ is a power series $F\in R[[X,Y]]$
with zero constant term
such that the ``addition'' rule on $\tau R[[\tau]]$ given by
$$g\oplus_Fh=F(g,h)$$
makes $\tau R[[\tau]]$ into an abelian group with identity~0.
In other words, for every $g,h$ we must have
$(f\oplus_Fg)\oplus_Fh=f\oplus_F(g\oplus_Fh)$ (associative law), $f\oplus_Fg=g\oplus_Ff$ (commutative
law), $f\oplus_F0=f$ (0 is identity), and for each $f\in\tau R[[\tau]]$ 
there exists $g\in\tau R[[\tau]]$ such that $f\oplus_F g=0$ (inverses).
Denote this group by ${\cal C}(F)$.  An equivalent and more widely
known definition is the following: a formal group law over $R$ is
a power series $F(X,Y)\in R[[X,Y]]$ such that \hfil\break
$$\eqalign{& ({\it i})\cr &({\it ii})\cr &({\it iii})\cr}\quad
\eqalign{& F(X,0)=X;\cr
&F(X,Y)=F(Y,X)\cr
&F(F(X,Y),Z)=F(X,F(Y,Z))\cr} \qquad
\eqalign{&{\rm (Additive\ Identity)}\cr &{\rm (Commutative\ Law)}\cr
&{\rm (Associative\ Law)}.\cr}\eqno(1.1)$$

The first property implies that $F$ has the form
$X+YH(X,Y)$.  By symmetry in $X$ and $Y$, it must therefore
be of the form
$$F(X,Y)=X+Y+XYG(X,Y), \qquad G\in R[[X,Y]].\eqno(1.2)$$

\proposition{1.1} 
{
Let $F$ be a power series in two variables with coefficients in $R$
such that $F(0,0)=0$.  The following are equivalent.
\hfil\break
(1) The three conditions in (1.1) hold;\hfil\break
(2) The binary operation on $\tau R[[\tau]]$ defined by $f\oplus_F g=
F(f,g)$ makes $\tau R[[\tau]]$ into an abelian group with identity 0;\hfil\break
(3) The binary operation on $\tau R[[\tau]]$ defined by $f\oplus_F g=
F(f,g)$ makes $\tau R[[\tau]]$ into an abelian semigroup with identity 0.
}

\proof{}   We will show $(1)\Rightarrow (2) \Rightarrow (3)
\Rightarrow (1)$. Assume (1) holds. Define a binary
operation on $\tau R[[\tau]]$ by $f\oplus_F g=F(f,g)$
for $f,g\in \tau R[[\tau]]$.  The three conditions immediately imply
$f\oplus_F 0=f$, $f\oplus_F g=g\oplus_F f$, and $(f\oplus_F g)\oplus_F h
=f\oplus_F(g\oplus_F h)$ for $f,g,h\in \tau R[[\tau]]$.  It remains
only to prove the existence of inverses.  For this, it suffices to prove
there is a power series $\iota\in\tau R[[\tau]]$ such that
$F(g,\iota\circ g)=0$ for all $g\in\tau R[[\tau]]$.
Let $\iota^{(1)}=-\tau$.
By (1.2)
$F(\tau,\iota^{(1)}) \equiv \tau-\tau\equiv0\bmod{\tau^2}$.
Now assume inductively that $\iota^{(N)}\in\tau R[[\tau]]$ satisfies
$F(\tau,\iota^{(N)}) \equiv0\bmod{\tau^{N+1}}$ and
$\iota^{(N)}\equiv \iota^{(N-1)}\bmod{\tau^{N}}$. Then there is $a\in R$
such that
$F(\tau,\iota^{(N)}) \equiv a\tau^{N+1}\bmod{\tau^{N+2}}$.
Let $\iota^{(N+1)}=\iota^{(N)}-a\tau^{N+1}$. By (1.2)
$$F(\iota^{(N)},-a\tau^{N+1})\equiv
\iota^{(N)}-a\tau^{N+1}=\iota^{(N+1)}
\bmod{\tau^{N+2}}.$$
Thus
$$\eqalign{F(\tau,\iota^{(N+1)})&\equiv
F(\tau,F(\iota^{(N)},-a\tau^{N+1}))
=F(F(\tau,\iota^{(N)}),-a\tau^{N+1})\cr
&\equiv F(\tau,\iota^{(N)})-a\tau^{N+1}\equiv 0\bmod{\tau^{N+2}}.\cr}$$
This completes the induction.  Let $\iota\in\tau R[[\tau]]$ be the
power series
such that $\iota\equiv \iota^{(N)} \bmod{\tau^{N+1}}$ for all $N$.  Then
$F(\tau,\iota(\tau))=0$, and hence $F(x,\iota(x))=0$ for all
$x\in\tau R[[\tau]]$.  This proves $(1)\Rightarrow (2)$.  It is obvious
that $(2) \Rightarrow (3)$.

Now assume (3) holds.  We will prove
condition ({\it iii}) of (1.1) holds; the other conditions in (1.1)
can be proved similarly.  Let $G(X,Y,Z)=F(F(X,Y),Z)-F(X,F(Y,Z))$.
We must show $G=0$.
By hypothesis, if $a,b,c$ are any positive integers then
$$G(\tau^a,\tau^b,\tau^c)=(\tau^a\oplus_F\tau^b)\oplus_F\tau^c-
\tau^a\oplus_F(\tau^b\oplus_F\tau^c)=0$$
as an element of $R[[\tau]]$.
We must show that every coefficient of $G$ is zero.  Write
$$G=\sum_{i,j,k\ge0} g_{ijk}X^iY^jZ^k.$$
Since the $N$th coefficient of $G(\tau^a,\tau^b,\tau^c)$ is zero
we have
$$\sum_{\{\,i,j,k\in\Z_{\ge0}\,|\,(a,b,c)\cdot(i,j,k)=N\,\}}g_{ijk}
=0\eqno(1.3)$$
for all positive integers $a,b,c,N$.
We need to show each $g_{ijk}=0$.  Suppose not.
Among all $i,j,k$ for which $g_{ijk}$ is nonzero,
consider those for which
$N_1=i+j+k$ is minimal.  Among all $i,j,k$ with
$g_{ijk}\ne0$ and $i+j+k=N_1$, consider those for which
$N_2=i+j$ is minimal.  Finally, among all $i,j,k$ with
$g_{ijk}\ne0$, $i+j+k=N_1$, and $i+j=N_2$ select the one
for which $N_3=i$ is minimal.  Call this triple
$(i_0,j_0,k_0)$; that is, $i_0+j_0+k_0=N_1$, $i_0+j_0=N_2$, $i_0=N_3$.
Choose integers $M_1,M_2,M_3$ such that
$$M_3\ge1,\quad M_2>M_3 N_3,\quad M_1> M_2N_2+M_3N_3.$$
Let
$$(a,b,c)=(M_1+M_2+M_3,M_1+M_2,M_1),\qquad N=M_1N_1+M_2N_2+M_3N_3.$$
We will obtain a contradiction by showing that
$$\sum_{\{\,i,j,k\in\Z_{\ge0}\,|\,(a,b,c)\cdot(i,j,k)=N\,\}}g_{ijk}
=g_{i_0,j_0,k_0}\ne0.\eqno(1.4)$$
Suppose $g_{ijk}\ne0$ and $(a,b,c)\cdot(i,j,k)=N$.  The equality
can be written
$$M_1(i+j+k)+M_2(i+j)+M_3i=N.\eqno(1.5)$$
Now $i+j+k\ge N_1$ by the minimality of $N_1$. Strict inequality
cannot hold, since otherwise
$$\eqalign{N&=M_1(i+j+k)+M_2(i+j)+M_3i\cr
&\ge M_1(N_1+1)
>M_1N_1+M_2N_2+M_3N_3
=N.\cr}$$
Thus $i+j+k=N_1$. By minimality of $N_2$ we know $i+j\ge N_2$.
Again strict inequality cannot hold, since otherwise
$$\eqalign{N&=M_1(i+j+k)+M_2(i+j)+M_3i\cr
&\ge M_1N_1+M_2(N_2+1)\cr
&>M_1N_1+M_2N_2+M_3N_3=N.\cr}$$
Thus $i+j=N_2$.  Now the equality (1.5) shows $i=N_3$.
This establishes (1.4) and completes the proof.\qed

The following proposition gives a general method to construct formal
group laws.  

\proposition {1.2} 
{
Let $G$ be an abelian group, $0_G$ its identity
element, and write its multiplication law additively.
Suppose there is a one-to-one map
$T:\tau R[[\tau]]\to G$ such that $T(0)=0_G$, and 
a power series $F\in R[[X,Y]]$ with zero constant term
such that
$$T(g)+T(h)=T(F(g,h))\eqno(1.6)$$
for all $g,h\in \tau R[[\tau]]$. Then 
$F$ defines a formal group law.
}

Some easy examples of the above proposition are: 
(1) $G=R[[\tau]]$ under addition, $T=$ inclusion,
$F(X,Y)=X+Y$ (called the {\bd additive group law}), and (2) $G=R[[\tau]]^\times$
under multiplication, $T(g)=1+g$, $F(X,Y)=X+Y+XY$ (called the {\bd
multiplicative group law}).  A less trivial example is the construction
of the group law associated to an elliptic curve, which will be given
in \S4.

\proof {of Proposition~1.2:}The hypothesis 
is that there is an injective map $T$ from $\tau R[[\tau]]$ into
an abelian group $G$ such that $T(0)=0_G$, and there is a power 
series $F(X,Y)$ with zero constant term such that 
$$T(g)+T(h)=T(F(g,h))$$
for all $g,h\in\tau R[[\tau]]$.  We need to show that $F$ gives
an abelian group law on $\tau R[[\tau]]$.  By Prop.~1.1,
it suffices to show $F$ makes $\tau R[[\tau]]$ into an
abelian semigroup with identity 0; that is, 
if $f,g,h\in\tau R[[\tau]]$ then
$$f\oplus_F(g\oplus_Fh)=(f\oplus_Fg)\oplus_Fh,\qquad f\oplus_F g=g\oplus_F f,\qquad f\oplus_F0=f.$$
Now $T(f\oplus_F(g\oplus_Fh))=T(f)+T(g\oplus_Fh)=T(f)+T(g)+T(h)$ and
similarly $T((f\oplus_Fg)\oplus_Fh)=T(f)+T(g)+T(h)$. This proves the
first identity, since $T$ is one-to-one.  The other
two identities are proved similarly.\qed

\titlea{2}{  Homomorphisms of formal group laws}

If $F$ is a formal group law then write ${\cal C}(F)$ for the 
group it determines.
That is, ${\cal C}(F)=\tau R[[\tau]]$ as a set, and the group 
law is given by $g\oplus_F h=F(g,h)$.
If $F,F'$ are two formal group laws then a {\bd homomorphism} from 
$F$ to $F'$ is defined as a power series 
$U(\tau)\in\tau R[[\tau]]$ with zero constant term
such that $g\mapsto U(g)$ defines a {\it homomorphism} from 
${\cal C}(F)$ into ${\cal C}(F')$.  Explicitly,
$$U\circ(x\oplus_F y)=(U\circ x)\oplus_{F'}(U\circ y)$$
for all $x,y\in\tau R[[\tau]]$. In terms of power series this
can be written
$$U(F(X,Y))={F'}(U(X),U(Y)).\eqno(2.1)$$
The reason that $U$ has zero constant term is that $U$ must
take $\tau R[[\tau]]$ into itself.
An example of a homomorphism from $F$ to itself
is the multiplication by $n$ map, denoted $[n]$ or $[n]_F$, 
which is defined by the rules:
$$\eqalign{[0]=0,\quad [1]=\tau,\quad 
&[n+1]\tau=[n]\tau\oplus_F \tau=F([n]\tau,\tau)
{\rm\ if\ } n>0,\cr
&[n]=\iota\circ[-n] {\rm\ if\ } n<0.\cr}\eqno(2.2)$$

Let $G_1,G_2$
be abelian groups, and let $T_i:\tau R[[\tau]]\to G_i$ ($i=1,2$)
be one-to-one maps such that $T_i(0)$ is the identity element of 
$G_i$.   Let $F_i$ be power series with zero constant term
such that
$$T_i(g)\oplus_{G_i}T_i(h)=T_i(g\oplus_{F_i}h),\qquad i=1,2,$$
where $\oplus_{G_i}$ denotes addition on the group $G_i$ and
$g\oplus_{F_i}h=F_i(g,h)$.
We showed that $F_i$ is a formal group law, and the above
equation simply states that $T_i$ is a group homomorphism
from ${\cal C}(F_i)$ into $G_i$.

\lemma{2.1}  
{
Let $G_i,T_i,F_i$, ${\cal C}(F_i)$ be as above.
Suppose there is a group homomorphism $\psi:G_1\to G_2$
and a power series $U$ with zero constant term such that 
$$\psi(T_1(g))=T_2(U(g))\eqno(2.3)$$
for all $g\in \tau R[[\tau]]$.  Then $U$ is a 
homomorphism between the formal group laws defined by $F_1$ and $F_2$.  
}

\proof{}
It suffices to show that $U$ is a homomorphism from ${\cal C}(F_1)$ to
${\cal C}(F_2)$.
By hypothesis there is a commutative diagram
$$\diagram{{\cal C}(F_1)&\lhook\joinrel\mathrel{\harr{T_1}{}}&G_1\cr
\varr{\displaystyle U}{}&&\varr{}{\displaystyle\psi}\cr
{\cal C}(F_2)&\lhook\joinrel\mathrel{\harr{T_2}{}}&G_2\cr}$$
Here $T_1,T_2,\psi$ are homomorphisms and $T_1,T_2$ are injective.
It follows by diagram chasing that $U$ is a homomorphism, as claimed.
\qed

As a special case, let $G_1=G_2=G$, $T_1=T_2=T$, $F_1=F_2=F$,
and $\psi(g)=ng$, where $n\in\Z$.  Then $U=[n]$, which was
defined by (2.2).  The power series for $[n]$ may either
be computed from the recursion (2.2) or from the formula (2.3),
which in this context reads
$$nT(g)=T([n](g))\qquad{\rm for\ }g\in\tau R[[\tau]].\eqno(2.4)$$
For the additive formal group law we have $T=$ inclusion of
$\tau R[[\tau]]$ into $R[[\tau]]$ and the formula reads
$ng=[n](g)$.  So in that case, 
$$[n](\tau)=n\tau \qquad{\rm (Additive\ Formal\ Group)}$$
For the multiplicative
formal group law we have $G= R[[\tau]]^\times$ and $T(g)=1+g$, so
the formula reads $(1+g)^n=1+[n](g)$.
In the special case where $n=p=$ the characteristic of $R$ with $p>0$
we have $(1+g)^p=1+g^p$, and therefore
$$[p](\tau)=\tau^p \qquad{\rm (Multiplicative\ Formal\ Group\ in\
Char.\ p).}$$

\titlea{3}{  Height}
If $R$ has characteristic $p$ then the {\bd height of a homomorphism}
$U$, written ${\rm ht}(U)$,
is the largest integer $h$ such that $U(\tau)=V(\tau^{p^h})$ for
some power series $V$, or $\infty$ if $U=0$.  
The {\bd height of the formal group law} is defined as
the height of the homomorphism $[p]$.  For 
the additive formal group law 
defined by $F(X,Y)=X+Y$ we have $[p](\tau)=p\tau=0$, 
so the height of $F$ is 
$\infty$.  For the multiplicative formal group law 
given by $F(X,Y)=X+Y+XY$ we have $[p](\tau)=\tau^p$,
therefore the multiplicative formal group law has height one.

\titlec{Example 3.1} Let $F=\sum f_{ij}X^iY^j$ be a formal
group law over an integral domain $R$ of characteristic~$p>0$.
Let $F^{(p)}=\sum f_{ij}^p X^iY^j$. We claim that $F^{(p)}$ is a 
formal group law, and $\phi=\tau^p$ is a homomorphism (evidently of
height~1) from $F$ to $F^{(p)}$.  For the first assertion,
replace $X,Y,Z$ by $X^{1/p}$, $Y^{1/p}$, $Z^{1/p}$ in the relation
(1.1) then take the $p$th power.  This yields the corresponding
relations for $F^{(p)}$.  For the second assertion,
note that 
$$F^{(p)}\bigl(\phi(X),\phi(Y)\bigr)=F(X,Y)^p
=\phi\bigl(F(X,Y)\bigr).$$
Observe that $\phi^k:F\to F^{(p^k)}$.
\qed

\proposition {3.2} {
Let $F_1,F_2$ be  formal group laws
over an integral domain $R$ of characteristic~$p$. Let
$U(\tau)=\sum u_i\tau^i$ be a homomorphism from 
$F_1$ to $F_2$ of height $k$.
Then the first nonzero coefficient of $U$ is $u_{p^k}$.
Moreover, there is a homomorphism $V:F_1^{(p^k)}\to F_2$ such that
$U=V\circ \phi^k$.
}

\proof{}  If $k=0$ then
$u_j\ne0$ for some $j$ which is prime to $p$, therefore
$U'(\tau)=\sum_m m u_m\tau^{m-1}$ is nonzero.  Differentiate
the equation $U(F_1(X,Y))=F_2(U(X),U(Y))$ with respect to $Y$
and then set $Y=0$.  We obtain
$$U'\bigl(F_1(X,0)\bigr){\partial F_1\over\partial Y}(X,0)
={\partial F_2\over\partial Y}\bigl(U(X),U(0)\bigr)\,U'(0).$$
Since $F_i(X,Y)=X+Y+XYG_i(X,Y)$ for $i=1,2$, this becomes
$$U'(X)\bigl(1+XG_1(X,0)\bigr)=\bigl(1+U(X)G_2(U(X),0)\bigr)\, u_1.$$
The left side is nonzero, therefore $u_1\ne0$.

Now let $k\ge1$ and set $q=p^k$.
By definition of height, there is a power series
$V(\tau)\in\tau R[[\tau]]$ such that $U(\tau)=V(\tau^q)$.
Now $V'$ is nonzero, since otherwise $V$ would be a function of 
$\tau^p$, so that $q$ could be replaced by $pq$.
We claim $V$ is a 
homomorphism from $F_1^{(q)}$ to $F_2$.  We have to show
$V\bigl(F_1^{(q)}(X,Y)\bigr)=F_2\bigl(V(X),V(Y)\bigr)$.
The left side is $V(F_1(X^{1/q},Y^{1/q})^q)
=U\bigl(F_1(X^{1/q},Y^{1/q})\bigr)$.
The right side is $F_2\bigl(U(X^{1/q}),U(Y^{1/q})\bigr)$.
These two are equal because $U$ is a homomorphism from $F_1$
to $F_2$. Since $V'\ne0$, $V$ has height zero. It follows from
the case $k=0$ that the first coefficient of $V$ is nonzero.
Thus the coefficient of $\tau^q$ in $U$ is nonzero.\qed

\proposition {3.3} {
Let $F,F',F''$ be formal group laws over
an integral domain $R$ of characteristic~$p$. In parts (a), (b),
(d) and (e) assume $p>0$. \hfil\break
(a) If $U:F\to F'$, and $V:F'\to F''$, then
${\rm ht}(V\circ U)={\rm ht}(V)+{\rm ht}(U)$.\hfil\break
(b) If there is a nonzero homomorphism
$U$ from $F$ to $F'$ then $F$ and $F'$ have the same height.\hfil\break
(c)  For $n\in \Z$, $[n]_F=n\tau+\tau^2(\cdots).$\hfil\break
(d) Every formal group $F$ over a ring of characteristic $p$ has 
height at least one.  \hfil\break
(e) If $n=ap^t$ with $(a,p)=1$ then ${\rm ht}([n]_F)=t\,\,{\rm ht}(F)$.
}

\proof{}  (a) Define the degree of a nonzero power series $\sum a_i\tau^i$
to be the smallest $i$ such that $a_i\ne0$. Prop.~3.2 asserts
that if $U$ is a nonzero homomorphism of formal group laws then
$\deg(U)=p^{{\rm ht}(U)}$.  
The degrees of power series multiply when they are composed, therefore
$p^{{\rm ht}(V\circ U)}=p^{{\rm ht}(V)}p^{{\rm ht}(U)}
=p^{{\rm ht}(V)+{\rm ht}(U)}$.
(b) Certainly $[p]_{F'}\circ U=U\circ[p]_F$, so 
$[p]_F$ and $[p]_{F'}$ have the same height by (a). (c) can easily be
shown by induction, using (2.2). (d) is immediate from (c) and
Prop.~3.2.
(e) ${\rm ht}([n]_F)={\rm ht}([a]_F)+t\,\,{\rm ht}([p]_F)$
by (a).  The height of $[a]_F$ is zero by (c),
and ht$([p]_F)={\rm ht}(F)$ by definition. 
\qed

If $F,F'$ are formal group laws over an integral domain $R$
and $U_1,U_2:F\to F'$, define $U_1\oplus_{F'}U_2=F'(U_1,U_2)$.  
$U_1\oplus_{F'}U_2$ is a homomorphism from $F$ to $F'$.
This composition rule makes Hom$(F,F')$ into an abelian group.
In particular, it is a $\Z$-module.  Suppose that $R$ has characteristic
$p>0$.  
We put a topology on ${\rm Hom}(F,{F'})$ by decreeing that $U$ and $V$
are close iff $U\ominus_{F'}V$ has a large height.  In other words,
the topology on $\Hom(F,{F'})$ is induced from the {\bd height metric}
$|U|=c^{{\rm ht}(U)}$, where $0<c<1$.\qed

\proposition {3.4}
{
  Let $F,F'$ be formal groups over an integral
domain $R$ of characteristic $p>0$.\hfil\break
(a) 
${\rm ht}(U_1\oplus_{F'}U_2)\ge{\rm inf}\{\,{\rm ht}(U_1),{\rm ht}(U_2)\,\}$.
If ${\rm ht}(U_1)<{\rm ht}(U_2)$ then ht$(U_1\oplus_{F'}U_2)={\rm ht}(U_1)$.
Hence, the height metric is nonarchimedean.\hfil\break
(b)  The map $\Z\times{\rm Hom}(F,{F'})\to{\rm Hom}(F,{F'})$
given by $(n,U)\mapsto [n]_{F'}\circ U$ is continuous with respect
to the $p$-adic metric on $\Z$ and the height metric
on ${\rm Hom}(F,{F'})$. Hence, Hom$(F,F')$ is naturally a $\Z_p$-module.
\hfil\break
(c) If ht$(F)<\infty$ then Hom$(F,F')$ is a faithful $\Z_p$-module.
}

\proof{} (a) Write ${F'}(X,Y)=X+Y+XYG'(X,Y)$. 
Then $U_1\oplus_{F'}U_2={F'}(U_1,U_2)=U_1+U_2+U_1U_2G'(U_1,U_2)$.
Part (a) is therefore true when the word ``degree''
is substituted for the word ``height''.  Since 
ht$(U_i)=\log_p(\deg(U_i))$, (a) follows.
(b)  We must show that if $n=m+ap^t$ with $t$ large
and if $U,V\in\Hom(F,{F'})$ are close then
$n\cdot U$ is close to $m\cdot V$.  But
$$n\cdot U\ominus_{F'}m\cdot V
=[n]_{F'}\circ(U\ominus_{F'}V)\oplus_{F'}[ap^t]_{F'}\circ V.$$
The height of $[n]_{F'}\circ(U\ominus_{F'} V)$ is $\ge{\rm ht}(U\ominus_{F'}V)$.
The height of $[ap^t]_{F'}\circ V$ is $\ge t$.  Both these heights
are large, so the height of the sum is large by (a). 
(c)  We must show that if $a\in\Z_p$ and $0\ne U\in{\rm Hom}(F,F')$ then
$a\cdot U=0$ iff $a=0$. Write $a=p^kb$, where $b\in\Z_p^\times$.
We have $a\cdot U=[p^k] \circ b\cdot U$. Certainly $b\cdot U\ne0$,
since $b$ is invertible, and $[p^k]$ is nonzero since
it has finite height. Thus $a\cdot U$ is the
composition of two nonzero formal power series over $R$, and since
$R$ is an integral domain, this composition is nonzero. \qed

It is a theorem of M.~Lazard ([3], [4])
that if $R$ is a separably closed field 
of characteristic $p$ then two formal group laws $F,F'$ 
defined over $R$
are isomorphic iff they have the same height; this gives a partial
converse to Prop.~3.3(b).  
We will see that the height of the formal group law 
associated to an elliptic curve $E$ defined over a field $R$ of 
characteristic $p$ is one or two according as $E$ is ordinary or 
supersingular.
Thus Lazard's Theorem implies that the formal group laws of any two
ordinary elliptic curves (or any two supersingular elliptic curves)
are isomorphic over the algebraic closure of $R$.  On the other hand,
the condition that
two elliptic curves over $R$ be isomorphic is much more restrictive
(the two curves must have the same $j$-invariant; see [7], p.~47-50)
This means that isomorphisms of formal group laws are 
far more abundant than isomorphisms of elliptic curves.

\titlea{4}{  Constructing the formal group law of an elliptic curve}

Let $E$ be an elliptic curve over a field $K$ determined by a nonsingular
Weierstrass equation 
$$W(X,Y,Z)=Y^2Z+a_1XYZ+a_3YZ^2-(X^3+a_2X^2Z+a_4XZ^2+a_6Z^3),\eqno(4.1)$$
$a_i\in K$.  Let $L$ be the quotient field of $K[[\tau]]$.
Since $K\subset L$, we can consider the points in $E(L)$.
Let $R$ be a subring of $K$ (possibly $R=K$) 
containing~1 and all the Weierstrass coefficients $a_i$.
We will construct a formal group law by embedding $\tau R[[\tau]]$
into $E(L)$ and ``stealing'' the group law from $E(L)$.

Consider points of the form $(t,-1,s)$ in $E(K)$.
Then $t$ can be regarded as the function $-X/Y\in K(E)$,
where $K(E)$ denotes the function field of $E$ over $K$,
and $t$ is a uniformizer at the identity $O=(0,1,0)$.  Also $s$ can be
regarded as the function $-Z/Y$, and $s$ has a triple zero
at $O$.  Let $\Omega$ be the
ring of functions in $K(E)$ which are defined at $O$
and $M$ the ideal of functions in $\Omega$ which vanish at
$O$.  Then $M$ is principal, generated by $t$, and $\Omega/M\cong K$
by the map $f+M\mapsto f(O)$.  
$\Omega$ has a metric induced by $M$, namely $|f|=c^{n}$,
where $0<c<1$ and $n$ is the largest integer such that $f\in M^n$.
The uniformizer $t$ determines an isometry
$\Psi:\Omega\to K[[\tau]]$ (where $K[[\tau]]$ has the $\tau$-adic
topology) as follows:
$f\mapsto \sum_{i=0}^\infty a_i\tau^i$ (where $a_i\in K$)
iff for each $N$, $f-\sum_{i=0}^N a_i t^i\in M^{N+1}$.
The image of $\Psi$ is dense in $K[[\tau]]$, since it contains
all polynomials.

Let $S(\tau)=\Psi(s)=\sum_{i=3}^\infty s_i\tau^i$.  We will prove below that
if $f\in\tau R[[\tau]]$ then $(f,-1,S(f))\in E(L)$,
so there is an embedding $T:\tau R[[\tau]]\to E(L)$ given by
$$T(f)=(f,-1,S(f)).\eqno(4.2)$$
The formal group law of $E$ will be the power series $F\in \tau
R[[\tau]]$ such that $T(g)+T(h)=T\bigl(F(g,h)\bigr)$.  All we need
to do is to prove this power series $F$ exists; it will
automatically be a formal group law because of Prop.~1.2.

By dividing through the Weierstrass equation by $Y^3$ we see that
$s$ and $t$ satisfy the equation 
$$s=t^3+a_1ts+a_2t^2s+a_3s^2+a_4ts^2+a_6s^3.\eqno(4.3)$$
The series $S$ can be computed by 
recursively substituting approximations
for $s$ into the right hand side of (4.3) 
and expanding to get improved approximations.
We start with the approximation $s=O(t^3)$ to obtain
$$\eqalign{s&
=t^3+a_1 t\,O(t^3)+a_2 t^2 O(t^3) +a_3(O(t^3))^2+a_4t(O(t^3))^2+a_6 (O(t^3))^3\cr
&=t^3+O(t^4).\cr}$$
On the next round substitute $t^3+O(t^4)$ for $s$ in the right side of the
equation to obtain $s=t^3+a_1t^4+O(t^5)$.
This procedure yields the general rule:
$$s_0=s_1=s_2=0,\qquad s_3=1,\qquad{\rm and\ if\ }n\ge4{\rm\ then}$$
$$s_n=a_1s_{n-1}+a_2 s_{n-2}+a_3\sum_{i+j=n} s_is_j
+a_4\sum_{i+j=n-1} s_is_j + a_6\sum_{i+j+k=n}s_is_js_k.\eqno(4.4)$$

\lemma {4.1}
{
 Let $W$ be the Weierstrass equation (4.1),
where $a_i\in R$ and $R$ is an integral domain.  Let $s_i\in R$
be defined by the recursion (4.4) and let
$S=\sum s_i\tau^i\in \tau R[[\tau]]$.  Then
$W(\tau,-1,S)=0$ in $R[[\tau]]$. If $f,g\in\tau R[[\tau]]$ and
$W(f,-1,g)=0$ then $g=S\circ f$.
}

\titled{ Remark.}  Since the Weierstrass equation is cubic
in the variable $Z$, it follows that for fixed $f\in\tau R[[\tau]]$,
the equation $W(f,-1,g)=0$ has three solutions
for $g$ in the algebraic closure of the quotient field of $R[[\tau]]$.
The lemma asserts that exactly one of these solutions lies in
$\tau R[[\tau]]$.
\proof{}  Let $K$ be the quotient ring of $R$ and 
let $E$ be the elliptic curve over $K$ with equation $W$.
Let $t=-X/Y$, $s=-Z/Y\in K(E)$, and $\Psi:\Omega\to K[[\tau]]$
be as described in the beginning of this section.
Then $\psi(t)=\tau$, $\Psi(s)=S$.  Now $W(t,-1,s)=0$, so
$$0=\Psi\left( W(t,-1,s)\right)=W(\tau,-1,S).$$
>From this it follows that $W(f,-1,S\circ f)=0$
for any $f\in \tau K[[\tau]]$.

Now suppose $f,g\in\tau R[[\tau]]$ and
$W(f,-1,g)=0$.  Let $h=S\circ f$. Then
$$\eqalign{0&=W(f,-1,h)-W(f,-1,g)\cr
&=(g-h)\,\left(-1+a_1f+a_2f^2+a_3(g+h) +a_4f(g+h)
+a_6(g^2+gh+h^2)\right). \cr}$$
Since $-1+a_1f+\cdots$ is a unit in $R[[\tau]]$, $g-h$
must be zero. \qed

The above lemma establishes that the map $T:\tau K[[\tau]]\to E(L)$
is well-defined, furthermore it is obviously one-to-one. 
Recall Prop.~1.2, which guarantees that if we can find
a power series $F$ in two variables with the properties
that $F(0,0)=0$ and $T(f)+T(g)=T(F(f,g))$ then $F$ will be a formal group law.
We now show such an $F$ can be found.
First we need to know addition formulas for points of the form 
$(t_1,-1,s_1)$.  Such formulas are provided below.

\proposition {4.2}
{
  Let $P_i=(t_i,-1,s_i)$ for $i=1,2$ be points
on the elliptic curve with Weierstrass equation (4.1). 
\hfil\break (a)\quad Suppose $t_1\ne0$ and let $m=s_1/t_1$.
If $1+a_2m+a_4m^2+a_6m^3\ne0$  then
$$-P_1=\left({-t_1\over 1-a_1t_1-a_3s_1},-1,{-s_1\over 1-a_1t_1-a_3s_1}
\right).\eqno(4.5)$$
(b)\quad Suppose $t_1\ne t_2$ and let $m=(s_1-s_2)/(t_1-t_2)$, 
$b=s_1-mt_1$, $A=1+a_2m+a_4m^2+a_6m^3$. 
If $A\ne0$ then
$$P_1+P_2=-(t_3,-1,mt_3+b),$$
$$t_3=-t_1-t_2-{a_1m+a_2b+a_3m^2+2a_4mb+3a_6m^2b\over A}.\eqno(4.6)$$
}

\proof{}
(b)\quad $P_1,P_2$ lie on the line
$mX-bY-Z=0$. Let $P_3$ be the third point of intersection of
this line with the elliptic curve. Write $P_3=(x_3,y_3,z_3)$.
If $y_3=0$ then $P_3=(1,0,m)$. From the Weierstrass equation
(4.1), $1+a_2m+a_4m^2+a_6m^3=0$, contrary to the hypothesis.
Thus $y_3\ne0$, and hence $P_3$
can be written $P_3=(t_3,-1,mt_3+b)$. Likewise $P_i=(t_i,-1,mt_i+b)$
for $i=1,2$.  When $(t,-1,mt+b)$ is substituted for $(X,Y,Z)$
in the Weierstrass equation, the result must be of the form
$A(t-t_1)(t-t_2)(t-t_3)$ with $A\ne0$.  Hence
$$\eqalign{-(mt+b)
&+a_1t(mt+b)+a_3(mt+b)^2+t^3+a_2t^2(mt+b)+a_4t(mt+b)^2\cr
&\qquad+a_6(mt+b)^3 = A(t-t_1)(t-t_2)(t-t_3).\cr}$$
The left side is of the form
$$\eqalign{(1+a_2m&+a_4m^2+a_6m^3)t^3+(a_1m+a_3m^2+a_2b+2a_4mb+3a_6m^2b)t^2\cr
&+(\cdots)t+(\cdots)\cr}$$
and the right side is of the form
$At^3-A(t_1+t_2+t_3)t^2+\cdots$. Now (b) follows immediately.

(a)\quad Let $P_2=(0,1,0)$,
$m=s_1/t_1$, $A=1+a_2m+a_4m^2+a_6m^3$.
Since $A\ne0$, (b) implies that
$P_1+(0,1,0)+(t_3,-1,mt_3)=(0,1,0)$, where $t_3=-t_1-(a_1m+a_3m^2)/A$.
Thus $-P_1=(t_3,-1,mt_3)$. Now
$$t_1^3A=t_1^3+a_2t_1^2s_1+a_4t_1s_1^2+a_6s_1^3=s_1-a_1t_1s_1-a_3s_1^2,$$
thus
$$\eqalign{t_3&=-t_1-{a_1m+a_3m^2\over
A}={-t_1(t_1^3A)-(a_1t_1^2s_1+a_3t_1s_1^2)\over t_1^3A}\cr
&={-t_1s_1\over s_1-a_1t_1s_1-a_3s_1^2}
={-t_1\over 1-a_1t_1-a_3s_1}.\cr}$$
\qed

\theorem {4.3}
{
  There is a power series $F(t_1,t_2)\in R[[X,Y]]$
with zero constant term such that for $f,g\in \tau R[[\tau]]$,
$$T(f)+T(g)=T(F(f,g)).\eqno(4.7)$$
Therefore $F$ is a formal group law.
}

\proof{}  
Consider Prop.~4.2, but treat $t_1,t_2$ as indeterminates
and substitute 
$S(t_1)$, $S(t_2)$ for $s_1,s_2$. In other words, we are working
over the field $L'=$ the quotient field of $R[[t_1,t_2]]$.
We need to show $t_3$ of equation (4.6) is a
power series in $t_1,t_2$. Let $M$ be the ideal of $R[[t_1,t_2]]$
generated by $t_1$ and $t_2$.  That is, $M$ is the
set of elements $\mu\in R[[t_1,t_2]]$ for which $\mu(0,0)=0$.
If $\mu\in M$ and $u$ is a unit of
$R$ then $u+\mu$ is a unit in $R[[t_1,t_2]]$.
Now
$$\eqalign{m
&={S(t_1)-S(t_2)\over t_1-t_2}=\sum_{i=3}^\infty {s_i(t_1^i-t_2^i)\over t_1-t_2}\cr
&=\sum_{i=3}^\infty s_i(t_1^{i-1}+t_1^{i-2}t_2+\cdots+t_1t_2^{i-1}+t_2^{i-1})\cr
}$$
so $m$ belongs to $M^2$.  Then $A=1+a_2 m+a_4m^2+a_6m^3$
is a unit in $R[[t_1,t_2]]$, since $A$ is the sum of a unit in $R$
and an element of $M$.  In particular, $A\ne0$, so Prop.~4.2(b)
applies. Also $b=S(t_1)-m t_1\in M^3$.  Now (4.6)
shows that $t_3\in M$.  Thus we can write $t_3=G(t_1,t_2)$, $G\in M$.
Certainly $t_3\ne 0$, because $G\equiv -t_1-t_2\bmod{M^2}$.
We have $(t_1,-1,S(t_1))+(t_2,-1,S(t_2))=-(t_3,-1,s_3)$ in $E(L')$,
where $s_3=mt_3+b\in M^3$.
By Prop.~4.2(a), the right side is
$$\left({-t_3\over 1-a_1 t_3-a_3s_3},-1,{-s_3\over 1-a_1t_3-a_3s_3}
\right).$$
Let 
$$F(t_1,t_2)={-t_3\over 1-a_1t_3-a_3s_3}\in M,\qquad
H(t_1,t_2)={-s_3\over 1-a_1t_3-a_3s_3}\in M^3.$$
If we substitute $t_1=f(\tau)$, $t_2=g(\tau)$ for $f,g\in\tau R[[\tau]]$
we get a homomorphism $R[[t_1,t_2]]\to R[[\tau]]$, which induces a
homomorphism $E(L')\to E(L)$. It follows that
$$(f,-1,S(f))+(g,-1,S(g))=(F(f,g),-1,H(f,g)).$$
By Lemma~4.1 $H(f,g)=S(F(f,g))$. This proves (4.7). The fact
that $F$ is a formal group law follows from Prop.~1.2.
\qed

The first few terms of $F$ are:
$$\eqalign{F(X,Y)&=X+Y-a_1XY-a_2(X^2Y+XY^2)\cr
&-(2a_3X^3Y+(3a_3-a_1a_2)X^2Y^2+2a_3XY^3)+\cdots\cr}$$

\titlea{5}{ Homomorphisms of formal group laws arising from isogenies}

Let $E,E'$ be two elliptic curves defined over the same field
$K$.  An {\bd algebraic map} from $E$ to $E'$ is a function
$\a:E(\cj K)\to E'(\cj K)$ such that for each $P\in E$ there exist
homogeneous polynomials $f_1,f_2,f_3$ of the same degree
and not all vanishing at $P$ such that for all but finitely many
$Q\in E(\cj K)$,
$$\a(Q)=(f_1(Q), f_2(Q),f_3(Q)).$$
An example of an algebraic map from $E$ to itself is the
{\bd translation by  $P$ map}
$\tau_P(Q)=P+Q$ for $P,Q\in E$.
The algebraic map is said to be {\bd defined} over a field $K$
if $E,E'$ are defined over $K$ and if all the coefficients of
$f_1, f_2, f_3$  can be chosen to belong to $K$.
It is a theorem ([7], p.~75) that every nonconstant algebraic
map from $E$ into $E'$ which takes the origin to the origin
is a group homomorphism.  Such an algebraic map is called an
{\bd isogeny}.
If $\tau:E\to E'$ and $-Q=\tau(0,1,0)\in E'$ then
$\tau_Q\circ\tau$ takes the origin of $E$ into the origin of
$E'$.
Thus every nonconstant algebraic map is
the composition of an isogeny with a translation.
Two curves $E,E'$ are called {\bd isogenous} over $K$
if there exists an isogeny defined over $K$
from $E$ into $E'$. The
{\bd endomorphism ring}  of $E$, written
$\End_K(E)$, is the set of isogenies over $K$ from $E$ to itself,
together with the constant zero map,
with the addition and multiplication
laws:
$$(\a+\b)(P)=\a(P)+\b(P),\qquad \a\b=\a\circ\beta.$$
Note that $\Z\subset\End_K(E)$. If $K$ is the finite
field with $q$ elements
then the {\bd Frobenius endomorphism} $\varphi_q$ is defined by
$\varphi_q(X,Y,Z)=(X^q,Y^q,Z^q)$.  Since $\varphi_q$
coincides with the Galois action, it commutes with any endomorphism
of $E$ which is defined over $K$. In particular, $\varphi_q$
commutes with $\Z$.

We claim that an isogeny of elliptic curves over $K$
gives rise to a homomorphism of the corresponding formal group laws 
over $K$.  Indeed, let
$$I(X,Y,Z)=(f_1(X,Y,Z),f_2(X,Y,Z),f_3(X,Y,Z))$$
be an isogeny between elliptic curves $E,E'$ over $K$. Here 
$f_1,f_2,f_3$ are homogeneous polynomials
of the same degree, say $d$, and 
$f_1,f_2,f_3$ do not simultaneously
vanish at the origin. Since the origin of $E$ is carried to the origin
of $E'$, $f_1$ and $f_3$ vanish at $O=(0,1,0)$ but $f_2(O)\ne0$.
Thus $f_1/Y^d\in M$ and $f_2/Y^d\in \Omega^\times$.
Now $f_1/Y^d=f_1(X/Y,1,Z/Y)=f_1(-t,1,-s)=(-1)^d f_1(t,-1,s)\in M$
and similarly $f_2/Y^d=(-1)^d f_2(t,-1,s)\in\Omega^\times$.  Thus
$$f_1(X,Y,Z)/f_2(X,Y,Z)=f_1(t,-1,s)/f_2(t,-1,s)\in M.$$ 
Let $U(\tau)=\sum_{i=1}^\infty u_i\tau^i$ 
denote the expansion of $f_1/f_2$
with respect to $t$.  Practically speaking,
$U$ can be obtained by expanding $s$ as a power
series $S$ and then computing 
$$f_1(\tau,-1,S(\tau))/f_2(\tau,-1,S(\tau))$$ 
in the ring
$K[[\tau]]$.  Note that $f_2(\tau,-1,S(\tau))$ is invertible since its
constant term is nonzero.

\proposition {5.1}
{
 Let $E,E',E''$ be elliptic curves over $K$
and let
$F,F',F''$ denote the associated formal group laws over $K$.
If $I:E\to E'$ is an isogeny then
the power series $U$ constructed above belongs to $\Hom(F,F')$.
The map $I\mapsto U$ is a one-to-one group homomorphism
from  
${\rm Isog}(E,E')\hookrightarrow \Hom(F,F')$. If $I':E'\to E''$
and $I'$ corresponds to $U'\in\Hom(F',F'')$ then
$I'\circ I$ corresponds to $U'\circ U\in\Hom(F,F'')$.
}

\proof{}  Let $L$ be the quotient field of $K[[\tau]]$.
Since $I$ is defined over $K$, it is a priori defined over $L$. 
The discussion above
shows that $I$ can be written in a neighborhood of the origin
as
$$I(X,Y,Z)=\left({f_1(t,-1,s)\over f_2(t,-1,s)},-1,{f_3(t,-1,s)\over
f_2(t,-1,s)}\right).$$
Let $T:\tau K[[\tau]]\to E(L)$ and 
$T':\tau K[[\tau]]\to E'(L)$ be the embeddings (4.2).
Substitute $(X,Y,Z)\to T(f)=(f,-1,S(f))\in E(L)$, where 
$f\in\tau K[[\tau]]$. Then $t=-X/Y$ changes to $f$ and $s=-Z/Y$
changes to $S\circ f$.
Thus 
$I(T(f))=(U(f),-1,V(f)),$
where 
$U(\tau)=f_1(\tau,-1,S(\tau))/f_2(\tau,-1,S(\tau))\in \tau K[[\tau]]$
and $V(\tau)=f_3(\tau,-1,S(\tau))/f_2(\tau,-1,S(\tau))\in \tau K[[\tau]]$.
By Lemma~4.1, $V=S'\circ U$, where $S'(t)$ is the power series
expansion for $-Z/Y$ in the curve $E'$.  Thus
$$I(T(f))=T'(U(f)).\eqno(5.1)$$
By Lemma~2.1,
this equation proves that $U$ is a homomorphism of formal group laws.

If $I_1,I_2\in{\rm Isog}(E,E')$, and if $U_1,U_2\in\Hom(F,F')$ are
the corresponding homomorphisms of formal group laws then on the
elliptic curve $E(L)$,
$$\eqalign{(I_1+I_2)&(\tau,-1,S\bigl(\tau)\bigr)\cr
&=I_1\bigl(\tau,-1,S(\tau)\bigr)+I_2\bigl(\tau,-1,S(\tau)\bigr)\cr
&=T'(U_1)+T'(U_2)\cr
&=T'(F'(U_1,U_2))\cr}\qquad
\eqalign{&\cr&{\rm by\ definition\ of\ }I_1+I_2\cr
&{\rm by\ (5.1)}\cr
&{\rm by\ (4.7).}\cr}$$
On the other hand, if $I_1+I_2$ corresponds to $U_3$ then
$$(I_1+I_2)\bigl(\tau,-1,S(\tau)\bigr)=T'(U_3).$$
Since $T'$ is one-to-one, $U_3=F'(U_1,U_2)=U_1\oplus_{F'}U_2$.
This shows that the map $I\mapsto U$ is a group homomorphism.

Finally, if $I:E\to E'$, $I':E'\to E''$ correspond to $U,U'$, 
respectively, then
since $U$ is the unique solution in $\tau K[[\tau]]$
to $I\circ T=T'\circ U$, 
$$I'\circ I\circ T
=I'\circ T'\circ U=T''\circ U'\circ U,$$
whence $I'\circ I$ corresponds to $U'\circ U$.
\qed
\titlec{ Example 5.2} Let $F$ be the formal group law over
$R$
associated to an elliptic curve $E$ with Weierstrass equation
(4.1), where the coefficients $a_i\in R$, and $R$ is an integral
domain. We will compute $[-1]_F$.
Let $g\in\tau R[[\tau]]$. By Proposition~4.2(a),
$$[-1]_E T(g)=[-1]_E(g,-1,S\circ g)
=\left({-g\over 1-a_1g-a_3 S\circ g},-1,
{-S\circ g\over 1-a_1g-a_3 S\circ g}\right)$$
The right side is $T\bigl(-g/(1-a_1g-a_3S\circ g)\bigr)$
by Lemma~4.1.  Now Lemma~2.1 implies
$$[-1]_F={-\tau\over 1-a_1\tau-a_3S}=-\tau\sum_{n=0}^\infty
(a_1\tau+a_3S)^n.$$
$\phantom{get square on right}$ \qed

An isogeny $I:E\to E'$ is called {\bd separable} 
if it  has the property: if $t'$ is a uniformizer at the origin of
$E'$ then $t'\circ I$ is a uniformizer at the origin of $E$.
This definition does not depend on the choice of uniformizer $t'$.
An isogeny which is not separable is called {\bd inseparable}. 
In characteristic zero, all isogenies are separable.
In characteristic $p$, the Frobenius is not separable, 
since it carries
uniformizers into $p$th powers of uniformizers.
It is a theorem ([7], II.2.12) that every isogeny can be factored as  
$\varphi_p^k$  from $E$ into $E^{(q)}$ $(q=p^k$)
composed with a separable isogeny from 
$E^{(q)}$ into $E'$. 

\lemma {5.3}
{
  Let $I$ be an isogeny from $E$ to $E'$ and
let $U(\tau)=\sum u_i\tau^i$ be the corresponding homomorphism 
between the formal group laws.  $I$ is separable iff $u_1\ne0$.
}

\proof{}  Let $t'$ be the function $-X/Y\in K(E')$.
$U$ is the power series expansion of $t'\circ I$ with respect to
the uniformizer $t=-X/Y\in K(E)$. Thus
$t'\circ I$ is not a uniformizer at the identity of $E$ iff
$t'\circ I\in M_{(0,1,0)}^2$ iff $u_1=0$.\qed

\titlec{ Example 5.4} Let $E$ be an 
elliptic curve whose Weierstrass coefficients $a_i$ belong to a field $K$
of characteristic~$p>0$, and let $F$ be its associated formal
group law.  Let $E^{(p)}$ be the elliptic curve with Weierstrass
coefficients $a_i^p$. Then the Frobenius map 
$\varphi_p:E\to E^{(p)}$ defined by $\varphi(X,Y,Z)=(X^p,Y^p,Z^p)$
corresponds to the
homomorphism of formal group laws $\phi=\tau^p:F\to F^{(p)}$.
\qed

\titlea{6}{  Height of an elliptic curve}

We begin this section with some facts about
elliptic curves over finite fields.
If $\a:E\to E'$ is an isogeny, define $\a^* K(E')=\{\,f\circ \a\,|\,
f\in K(E')\,\}$; this is a subfield of $K(E)$.
The {\bd degree} of an isogeny $\a:E\to E'$
is the index of $\a^* K(E')$ in $K(E)$.
This number is finite because both fields have transcendence
degree~1 and $\a$ is a nonconstant map.
If $K$ has characteristic $p$
then the Frobenius isogeny $\varphi_p(X,Y,Z)=(X^p,Y^p,Z^p)$
from $E$ into $E^{(p)}$ has degree $p$. Here
$E^{(p)}$ is the curve whose Weierstrass equation
is obtained from that of $E$ by raising the coefficients to the
$p$th power.

Every isogeny $\a:E\to E'$ has a {\bd dual isogeny} $\hat\a:E'\to E$.
The dual isogeny is characterized by the property that
$\a\circ\hat\a=[\deg(\a)]_{E'}$ and $\hat\a\circ\a=[\deg(\a)]_E$,
where $[n]_E$ denotes multiplication by $n$.  If $E=E'$, then
there is an integer $a(\a)$, called the {\bd trace of $\a$},
such that $\a+\hat\a=[a(\a)]_E$.  The endomorphism $\a$
satisfies the quadratic equation
$$\a^2-[a(\a)]\a+[\deg(\a)]=0\qquad{\rm in\ }\End(E).$$
In particular, if $K$ has $q$ elements then there is $t\in\Z$
such that
$$\varphi_q^2-[t]\varphi_q+[q]=0.$$
The integer $t$ is called the {\bd trace of Frobenius}.
It is well known ([7], Ch.~5) that $|t|\le 2\sqrt q$ and the cardinality
of $E(K)$ is $q+1-t$.

The height of a formal group law was defined in \S 3. 
Naturally, the height of an elliptic curve is defined to be
the height  of the associated formal group law.

\proposition {6.1}
{
  An elliptic curve over a field of
characteristic~$p$, where $p>0$, has height one or two.
}

\proof{}  Let $\varphi_p:E\to E^{(p)}$ be the $p$th power
Frobenius and $\hat\varphi_p:E^{(p)}\to E$ its dual.
Let $F$ be the formal group law associated to $E$, and
let $V(\tau)=\sum v_i\tau^i:F^{(p)}\to F$ 
be the homomorphism of formal group laws associated
to $\hat\varphi_p$.  Then $[p]_F=V(\tau^p)$. If
$\hat\varphi_p$ is separable then $v_1\ne0$, so $E$ has
height one. If $\hat\varphi_p$ is inseparable, it
can be written as a composition of a power of $\varphi_p$
and a separable isogeny ([7], Corollary II.2.12). 
Since the degree of $\hat\varphi_p$
equals the degree of $\varphi_p$, only one power of $\varphi_p$
can occur in this decomposition. Thus
$\hat\varphi_p=\a\circ\varphi_p$ with $\a$ an isomorphism.
Let $A=\sum a_i\tau^i$ be the power series corresponding to $\a$
and let $A'$ be the power series corresponding to $\a^{-1}$.
Then $[p]_E=A(\tau^{p^2})=a_1\tau^{p^2}+\cdots$, and $a_1\ne0$
because $A\circ A'(\tau)=\tau$.
In this case $E$ has height two.\qed

An elliptic curve in characteristic~$p$
of height one is called {\bd ordinary}.
An elliptic curve in characteristic~$p$
of height~2 is called {\bd supersingular}. The next lemma
gives another characterization of supersingular and ordinary curves
when the underlying field is finite. 

\proposition {6.2}
{
 An elliptic curve $E$ over a finite field $K$ with
$q=p^n$ elements is supersingular iff $p$ divides the trace of 
Frobenius iff $|E(K)|\equiv1\bmod p$. If $E$ is supersingular and
$n$ is even then $|E(K)|=q+1+m\sqrt q$, $m\in\{\,-2,-1,0,1,2\,\}$.
If $E$ is supersingular, $n$ is odd, and $p\ge 5$,  then
$|E(K)|=q+1$. If $E$ is supersingular, $n$ is odd, and $p\le 3$
then $|E(K)|=q+1+m\sqrt{pq}$, where $m\in\{\,-1,0,1\,\}$.
}

For a more precise statement about which values of $|E(K)|$
can occur, the reader may consult [8], Theorem~4.1.

\proof{} As above, let $F$ be the formal group law corresponding to $E$
and $V:F^{(p)}\to F$ the homomorphism of formal group laws corresponding
to $\hat\varphi_p$. In other words, $V$ is defined by 
$[p]_F=V(\tau^p)$.  Recall that $E^{(p)}$ denotes the elliptic
curve whose Weierstrass equation is obtained by
taking the $p$th powers of the Weierstrass coefficients for
$E$, and we use similar notation for isogenies.  Now
$\hat\varphi^{(p^k)}:E^{(p^{k+1})}\to E^{(p^k)}$ 
is the dual of the map $\varphi_p:E^{(p^k)}\to E^{(p^{k+1})}$, so 
$$\hat\varphi_p\circ \hat\varphi_p^{(p)}\circ\cdots\circ 
\hat\varphi_p^{(p^{n-1})}$$
is the dual of $\varphi_p^n$.  The corresponding formal group law 
homomorphism is
$$N(V)= V\circ V^{(p)}\circ\cdots\circ V^{(p^{n-1})}.$$
Let $t$ be the trace of Frobenius, so that
$|E(K)|=q+1-t$.  Since $[t]_E$ is the sum of
$\varphi_p^n$ and its dual in $\End(E)$, it follows that
$$[t]_F=N(V)\oplus_F \tau^{p^n}=F(N(V),\tau^{p^n}).$$
If $E$ is supersingular then $V$ has height one, so
${\rm N}(V)$ has height $n$. In that case, $[t]_F$ has
height at least $n$, so $[t^2]_F$ has height at least $2n$. 
Since the height of $F$ is two in this case,
Prop.~3.3(e) implies $t^2$ is divisible by $p^n$.
Since $|t|\le 2\sqrt q$ and $q|t^2$, we deduce that
$t^2\in\{\,0,q,2q,3q,4q\,\}$.  Since $t\in\Z$, we find
$t\in\{\,0,\pm q^{1/2},\pm2q^{1/2}\,\}$ if $n$ is even; 
$t=0$ if $n$ is odd
and $p>3$; $t\in\{\,0,\pm\sqrt{2q}\,\}$ if $n$ is odd and $p=2$,
$t\in\{\,0,\pm\sqrt{3q}\,\}$ if $n$ is odd and $p=3$.
Since $|E(K)|=q+1-t$,
the cardinality of $E(K)$ must be of the form stated.

Next suppose $E$ is ordinary. Then N$(V)$ has height zero,
so $[t]_F$ has height zero.  In that case Prop.~3.3(e) implies
$t$ is prime to $p$.\qed

\proposition {6.3}
{
 If $E$ is an ordinary
elliptic curve defined over a 
field $K$ of cardinality $p^n$ 
and $F$ is its associated formal group law then
the trace of the Frobenius endomorphism is equal mod~$p$
to the norm
from $K$ to $\F_p$ of the first nonzero coefficient of $[p]_F$.
}

\proof{}  Let $|K|=p^n=q$.  The homomorphism of $F$ associated to 
$\varphi_q^2+[-t]_E\varphi_q+[q]_E$
is zero, thus each of its coefficients is zero.
Now $\varphi_q$ corresponds to the power series $\tau^q$, and
$[-t]_E$ corresponds to a power series of the form 
$-t\tau+\tau^2(\cdots)$, therefore 
$\varphi_q^2+[-t]_E\circ\varphi_q$
corresponds to $F(\tau^{q^2},-t\tau^q+\tau^{2q}(\cdots))$,
which is of the form $-t\tau^q+\tau^{2q}(\cdots)$. Finally,
we evaluate $[q]_F$. Let $\phi=\tau^p$. Since $\phi\circ V=V^{(p)}\circ\phi$,
$$[q]_F=(V\circ\phi)^n=V\circ V^{(p)}\circ\cdots\circ 
V^{(p^{n-1})}\circ\phi^n=({\rm N}_{K/{\bf F}_p}(v)\tau+(\cdots)\tau^2)\circ\tau^q,$$
so $[q]_F={\rm N}_{K/{\bf F}_p}(v)\tau^q+(\tau^{2q})(\cdots)$.  Thus
$$0=F\left(-t\tau^q+\tau^{2q}(\cdots), {\rm N}_{K/{\bf F}_p}(v) \tau^q
+ \tau^{2q}(\cdots)\right)
=(-t+{\rm N}_{K/{\bf F}_p}(v))\tau^q+\tau^{2q}(\cdots).$$
\qed

\titlea{7}{  Some theorems of Couveignes}

Let $R$ be an integral domain of characteristic~$p$.
Let $\F_p\subset\R$ be the field
with $p$ elements if $p$ is prime, and $\F_p=\Z$ if $p=0$.
Let
$$F=\sum_{i,j}f_{ij}X^iY^j,\qquad F'=\sum_{i,j}f'_{ij}X^iY^j$$
be two formal group laws over $R$, and let $U(\tau)=\sum_{i=1}^\infty
u_i\tau^i\in\tau R[[\tau]]$ be a 
homomorphism from $F$ to $F'$.
Couveignes proved with an elementary argument in his PhD thesis
that the coefficients $u_i$ satisfy some
simple relations over $R$.  He used these relations to compute the
orders of elliptic curves over finite fields of small characteristic
(see [2] and [6]).
In [1] it is shown that Couveignes' method is closely related to the
modified Schoof algorithm which was developed by Atkins and Elkies;
see [5] and its bibliography.
In this section we state and prove Couveignes' theorems.  In the next
section we prove related results which are used in [1].

\theorem {7.1}
{
  Let $i$ be a positive integer which is
not a power of~$p$.  If $p=0$ assume $i\choose m$ is a unit in $R$
for some $1\le m<i$.
There is a polynomial $C_i$ in several variables
with coefficients in $\F_p$
such that for each $F,F',U$ as above we have
$$u_i=C_i(u_j,f_{k\ell},f'_{k\ell}\,|\,1\le j<i, 1\le k+\ell\le i\,).$$
}

\proof{}  Let $A$ be transcendental and work in the integral domain
$R[A]$.  Since $U$ is a homomorphism,
$$U(F(\tau,A\tau))=F'(U(\tau),U(A\tau)).$$
By (1.2) there are power series $G,G'\in R[[X,Y]]$
such that $F(X,Y)=X+Y+XYG(X,Y)$ and $F'(X,Y)=X+Y+XYG'(X,Y)$.
Therefore
$$\eqalign{\sum u_j&(\tau+A\tau+A\tau^2 G(\tau,A\tau))^j=\cr
&\sum u_j\tau^j+\sum u_j(A\tau)^j+U(\tau)U(A\tau)G'(U(\tau),U(A\tau)).\cr}$$
This can be rewritten
$$\eqalign{0=&\sum u_j\tau^j\{(1+A+A\tau G(\tau,A\tau))^j-(1+A^j)\}\cr
&- A\tau^2(\sum_{j=0}^\infty u_{j+1}\tau^j)
(\sum_{j=0}^\infty u_{j+1}(A\tau)^j)
G'(\sum_{j=1}^\infty u_j\tau^j,\sum_{j=1}^\infty u_j(A\tau)^j).
\cr}$$
The coefficient of $\tau^i$ is 
of the form $u_i\{(1+A)^i-(1+A^i)\}+M_i$, where $M_i$ is a 
polynomial in $A,u_1,u_2,\ldots,u_{i-1}$ and in some of the 
coefficients of
$G,G'$.  This gives the relation
$$u_i\{(1+A)^i-(1+A^i)\}-M_i =0.$$
The hypothesis that $i$ is not a power of $p$ implies
$(1+A)^i\ne 1+A^i$.  If $p=0$ choose $m$ such that $i \choose m$ is
a unit in $R$, and if $p>0$ let
$m$ be a positive integer such that the coefficient of $A^m$ is nonzero
in the polynomial $(1+A)^i-(1+A^i)$.  In characteristic~$p$
this coefficient is a unit in $R$ because it is a nonzero
element of the prime field $\F_p$.
Since $A$ is transcendental,
the coefficient of $A^m$ in our relation must be identically
zero. This coefficient gives our desired formula
for $u_i$ in terms of the $u_j$ and the coefficients of $F$ and $F'$.
\qed

The next theorem accounts for the $u_i$ when $i$ is a power of $p$.
It was proved by Couveignes for formal group laws associated
to ordinary elliptic curves, but his argument generalizes
easily to formal group laws of any height.  

\theorem {7.2}
{
 Let $i$ be a power of a prime $p$ and let 
$h>0$.  There is a 
polynomial $C_i$ in several variables with coefficients in
$\F_p$ such that:
if $F=\sum f_{k\ell}X^kY^\ell$ and ${F'}=\sum f'_{j\ell}X^jY^\ell$
are  formal group laws of height $h$ over a domain $R$ of characteristic
$p$ and $U=\sum u_j\tau^j:F\to {F'}$ a homomorphism
then 
$$ v'_1 u_i^q-v_1^i u_i
=C_i(u_j,f_{k\ell},f'_{k\ell}\,|\,j<i, k+\ell\le qi\,)$$
where $q=p^h$ and
$v_1,v_1'$ are the first nonzero coefficients of the power series
$[p]_F$, $[p]_{F'}$, respectively.
}

\proof{}  By Prop.~3.2 we can write 
$[p]_F(\tau)=V\circ\phi^h(\tau)=V(\tau^q)$,
where $V(\tau)=\sum v_j\tau^j$ 
is a homomorphism of height zero from $F^{(q)}$ to ${F'}$.
It is easy to show by induction on $n$ that for $n>0$ 
the $j$th coefficient of $[n]_F$ is a polynomial in the $f_{k\ell}$
with $k+\ell\le j$.  Since $v_j$ is the $jq$th coefficient
of $[p]_F$, $v_j$ is a polynomial in the $f_{k\ell}$ with 
$k+\ell\le jq$. Similarly $[p]_{F'}=V'\circ\phi^h$, 
$V'(\tau)=\sum v_j'\tau^j$,
and $v'_j$ is a polynomial in the $f'_{k\ell}$ with $k+\ell\le jq$.
Since  $[p]_{F'}\circ U=U\circ[p]_F$,
$$V'(U(\tau)^q)=U\bigl(V(\tau^q)\bigr).$$
Let $\s=\tau^q$. The left side is
$$v'_1(\sum_{j=1}^\infty u_j^q\sigma^{j})+
v'_2(\sum_{j=1}^\infty u_j^q\sigma^{j})^2+\cdots,$$
and the coefficient of $\sigma^{i}$ is of the form
$v'_1u_i^q$ plus terms involving $u_j$ for $j<i$ and $v'_j$ for $j\le i$.
The right side is
$$u_1(\sum_j v_j\sigma^{j})+u_2(\sum_j v_j\sigma^{j})^2+\cdots+
u_i(\sum_j v_j\sigma^{j})^i+\cdots.$$
This time the coefficient of $\sigma^{i}$ is of the form
$u_i (v_1)^i$ plus terms involving $u_j$ for $j<i$ and 
$v_j$ for $j\le i$.
By equating the two sides we get $v_1' u_i^q-v_1^i u_i$
equals a polynomial in the $u_j$ for $1\le j<i$
and the $v_j, v_j'$ for $1\le j\le i$.\qed

\titlea{8}{ Further results relating to  Couveignes' theorems}

Fix the following notation throughout this section.
Let $R$ be an integral domain of characteristic~$p>0$,
$F$ and ${F'}$ formal group laws of height $h$ over $R$, and $q=p^h$.
Let $C_1,C_2,\ldots$
denote Couveignes' relations given in \S 7 evaluated at the
coefficients of $F,{F'}$ but leaving the $u_i$ as indeterminates; 
thus $C_i\in R[X_1,\ldots,X_i]$
and $C_i=X_i+$ a certain polynomial 
in $X_1,\ldots,X_{i-1}$ if $i$ is not a power of $p$;
$C_i=v'_1X_i^q-v_1^iX_i+$ a certain polynomial in 
$X_1,\ldots,X_{i-1}$ if $i$ is a  power of $p$.
Here the $v_i$  and $v_i'$ lie in $R$, since they are 
polynomials in the coefficients of $F$ and $F'$, respectively.
Couveignes' theorems assert that if $\sum u_i\tau^i\in\Hom(F,{F'})$
then $C_i(u_1,\ldots,u_i)=0$ for all $i$.
Let $K$ denote the separable algebraic closure of the quotient
field of $R$.  

\lemma {8.1}
{
There are exactly $q^n$ solutions $(u_1,\ldots,u_{p^n-1})$
with $u_i\in K$ to the first $p^n-1$ of Couveignes' relations.
}

\proof{}  
For each solution $(w_1,\ldots,w_{i-1})$
to the first $i-1$ of Couveignes' equations over $K$ there are 
$q$ values or~1 value of $w_i$ such that $(w_1,\ldots,w_i)$ is a
solution to the $i$th relation, 
according as $i$ is or is not a power of $p$. 
(To see that the $q$ solutions for $w_i$
are distinct when $i$ is a power of $p$, note that the derivative
with respect to $X_i$ of $C_i$ is $v_1^i$, which is nonzero.)
The lemma now follows easily by induction on $n$.\qed

\theorem {8.2}
{
 If $u_1,u_2,\ldots$
is a solution to Couveignes' relations then 
$\sum u_i\tau^i\in\Hom(F,{F'})$.
}

\proof{}  Without loss of generality we can replace $R$ by $K$.
In Chapter III, \S 2 of [3] it is shown that $\Hom(F,{F'})$ is free over
$\Z_p$ of rank $h^2$ and $p^n\Hom(F,{F'})$ is the set of homomorphisms
with height $\ge nh$. (In fact, it is shown that
$\Hom(F,{F'})$ is the maximal order of
a central division algebra over $\Q_p$ of rank $h^2$ and
invariant $1/h$, but we do not need this here.)
It follows that a complete set of $\Z_p$-module generators
$U_1,\ldots,U_{h^2}$ can be found such that the height of each
generator is less than $h$, and if $\sum c_iU_i$ has height
$\ge nh$ for some $c_i\in\Z_p$ then each $c_i$ is divisible by $p^n$.
If $U,{U'}\in\Hom(F,{F'})$ and $U\equiv {U'}\bmod{\deg\,\,q^n}$ 
(meaning that the $i$th coefficient of $U$ and $U'$ coincide
for all $i\le q^n$) then
$$0={F'}({U'},[-1]_{F'}\circ {U'})\equiv {F'}(U,[-1]_{F'}\circ {U'})=U\ominus_{F'} {U'}
\bmod{\deg\,\, q^{n}},$$
so $U\ominus_{F'} {U'}$ has height $\ge nh$, 
and it is therefore divisible by $p^n$.  
Thus $\sum c_i U_i\equiv\sum c_i' U_i\bmod\deg\,\, q^{n}$
$(c_i,c_i'\in\Z_p$) implies $c_i\equiv c_i'\bmod p^n$.  This shows that
the number of distinct elements $\sum_{i=1}^{q^{n}-1} u_i\tau^i$
which are truncations of power series in $\Hom(F,{F'})$ is
the cardinality of $(\Z/p^n\Z)^{h^2}$, which is
$q^{nh}$.  Each truncation gives rise to a solution $(u_1,\ldots,
u_{q^{n}-1})$ of the first $q^{n}-1$ of Couveignes' relations.
Since this coincides with the total number of solutions,
each solution of Couveignes' relation arises from $\Hom(F,{F'})$.\qed

\corollary {8.3}
{
 If $h=1$
and if $\Hom(F,{F'})$ contains a homomorphism (with coefficients in $R$)
of height $k$
then all the solutions $(v_1,v_2,\ldots)$
in $K$ to Couveignes' relations for which $v_i=0$ for
$i<p^k$ actually lie in $R$.
}

\proof{}  Let $U$ be the homomorphism of height $k$ and
$\Z_p\cdot U=\{\,c\cdot U\,|\,c\in\Z_p\,\}$. As mentioned in
the previous proof, Hom$(F,{F'})\cong\Z_p$, and it is
generated by a homomorphism $U_0$ of height zero. 
Find $a\in\Z_p$ such that $U=a\cdot U_0$. Since ht$(a\cdot U_0)=
v_p(a)$, $v_p(a)=k$.  Thus $\Z_p\cdot U=\Z_pa\cdot U_0
=p^k\Z_p\cdot U_0$.
Since $U$ is defined over $R$, so is $c\cdot U$ for each $c\in\Z_p$.
Thus every element of $p^k\Z_p\cdot U_0$ has coefficients in $R$.
The coefficients of such elements are precisely the solutions
$(v_1,v_2,\ldots)$
to Couveignes' relations which have $v_i=0$ for all $i<p^k-1$.
\qed

\begref{References}{8.}

\refno {1.} A.~W.~Bluher, Relations between certain power sums of elliptic
modular forms in characteristic two, {\sl J.\ Number Theory} {\bf 70},
127-145 1998
\refno {2.} J.~M.~Couveignes, Quelques calculs en theorie des nombres,
Ph.D.~thesis, Bordeaux, 1995
\refno {3.} A.~ Frohlich, {\it Formal Groups}, 
{\sl Lect.\ Notes in Math.\ }{\bf74}, Springer-Verlag, 1968
\refno {4.} M.~Hazewinkel, {\it Formal Groups and Applications},
Academic Press, New York, 1978
\refno {5.} R.~Lercier and F.~Morain, Counting the number of points on
elliptic curves over $F_{p^n}$ using Couveignes' algorithm,
Research report LIX/RR/95/09, Ecole Polytechnique-LIX,
September 1995
\refno {6.} R.~Lercier and F.~Morain, Counting the number of points on
elliptic curves over finite fields: strategies and performances,
{\sl Advances in Cryptology -- EUROCRYPT '95} {\sl Lect.\ Notes in
Computer Science }{\bf921}, Springer, 1995, 79-94
\refno {7.} J. H. Silverman, {\it The Arithmetic of Elliptic Curves},
Springer-Verlag, New York, 1986
\refno {8.} W. C. Waterhouse, Abelian varieties over finite fields,
{\sl Ann.\ Scient.\ \'Ec.\ Norm.\ Sup.} {\bf 2} 1969, 521-560
\endref
\end